\documentstyle[11pt,leqno]{article}

\input amssym.def
\input amssym.tex
\font\Bbb=msbm10 
scaled\magstep1\font\scriptBbb=msbm10\font\scriptscriptBbb=msbm7\newfam
\Bbbfam\textfont\Bbbfam=\Bbb\scriptfont\Bbbfam=\scriptBbb
\scriptscriptfont\Bbbfam=\scriptscriptBbb

\def\CC{{\fam=\Bbbfam C}}

\def\HH{{\fam=\Bbbfam H}}

\def\NN{{\fam=\Bbbfam N}}

\def\RR{{\fam=\Bbbfam R}}

\def\TT{{\fam=\Bbbfam T}}

\newtheorem{prop}[subsection]{Proposition}

\newtheorem{remark}[subsection]{Remark}
\newtheorem{example}[subsection]{Example}
\newtheorem{theorem}[subsection]{Theorem}
\newtheorem{cor}[subsection]{Corollary}
\newtheorem{basic lemma}[subsection]{Basic Lemma}
\newtheorem{lemma}[subsection]{Lemma}

\def\lan{\langle}
\def\ran{\rangle}
\def\nat{\sharp}
\def\noi{\noindent}
\def\o{\overline}
\def\de{\partial}
\def\RE{\Re\mbox{\rm e}}
\def\IM{\Im\mbox{\rm m}}

\def\al{\alpha}
\def\eps{\varepsilon}
\def\la{\lambda}
\def\om{\omega}
\def\Om{\Omega}
\def\ph{\varphi}
\def\vr{\varrho}
\def\Lie{{\rm Lie}}
\def\TT{{\fam=\Bbbfam T}}

\def\h{{\frak h}}
\def\sp{\frak{sp}}
\def\u{\frak u}

\def\A{{\cal A}}
\def\B{{\cal B}}
\def\C{{\cal C}}
\def\cS{{\cal S}}
\def\D{{\cal D}}
\def\E{{\cal E}}
\def\H{{\cal H}}
\def\Aut{{\rm Aut}\,}
\def\id{{\rm id}}
\def\Sp{{\rm Sp}}
\def\sign{{\rm sign}}
\def\supp{{\rm supp}}
\def\trans{\,{}^t\!}
\def\tr{{\rm tr}}
\def\HS{{\cal H\cal S}}

\def\be{\begin{enumerate}} 
\def\ee{\end{enumerate}}

\begin{document}

\title{Analysis of invariant PDO's on the Heisenberg group\\ 
      \large ICMS-Instructional Conference, Edinburgh 7.-13.4.1999}
\author{Detlef~M\"uller}
\date{}
\maketitle

\section{Analysis on the Heisenberg group: Basic Facts}
\vskip-0.3cm\hskip1cm (see e.g.\cite{Fo}, \cite{T})

\subsection{The Heisenberg group and its automorphisms}
 
The {\it Heisenberg group } $\HH_n$ is $\RR^n\times \RR^n\times \RR$, endowed with 
the product 
\[
(x,y,u)\cdot (x',y',u'):=(x+x',y+y',u+u'+\frac 1 2 (x\cdot y'-y\cdot x')).
\]
Observe:
\begin{itemize}
\item 0 is the neutral element of $\HH_n$
\item $(x,y,u)^{-1}=(-x,-y,-u)$
\item Writing $z=(x,y)\in \RR^{2n}$, and regarding $z$ as a column vector, we may
 regard 
$\HH_n$ also as $\RR^{2n}\times \RR$, with product
\end{itemize} 
\begin{eqnarray}%1.1
(z,u)\cdot (z',u')=(z+z',u+u'+\frac 1 2 \lan z,z'\ran),
\end{eqnarray}
where $\lan\ ,\ \ran$ denotes the {\it canonical symplectic form}
\[
\lan z,w\ran :=\ \trans z\cdot J\cdot w, \qquad J=\pmatrix{ 0 & I_n\cr
        -I_n & 0}
\]
on $\RR^{2n}$.
\bigskip

\noi{\bf Exercise:} $\HH_n$ is isomorphic to the group of upper triangular matrices

\[
\pmatrix{ 1 & p_1 & \dots & p_n & t\cr
        & 1 & & 0 & q_1\cr
        &&\ddots &&\vdots\cr
        &&&1&q_n\cr
        &&&&1} , \qquad p_j,q_j,t\in \RR.
\]
If $\om$ is any symplectic bilinear form on a finite dimensional vector space 
$V$ over  $\RR$ or $\CC$, 
denote  by 
\begin{eqnarray}
Sp(\om):=\{T\in L(V,V): \ \om(Tz,Tw)=\om (z,w)\ \forall z,w\in V\}
\end{eqnarray}
the corresponding {\it symplectic group}. If $\om=\lan\ ,\ \ran$, we also 
write $Sp(n,\RR)$ resp. 
$Sp(n,\CC)$ for these groups. Notice:
\[
T\in Sp(n,\RR) \ \Longleftrightarrow \ ^tT\cdot J \cdot T=J
\]
If $t\mapsto T(t)$ is a smooth curve in $Sp(\om)$ with $T(0)=I$, one finds  
from (1.2) that $S:= \frac{dT}{dt} (0)$ 
satisfies
 
\begin{equation}%1.3
\om(Sz,w)+\om(z,Sw)=0,
\end{equation}
i.e. $S$ is skew symmetric w.r. to $\om$.

This shows that the Lie algebra $\sp(\om)$ of $Sp(\om)$ consists of all linear
 endomorphisms $S$ of $V$ satisfying (1.3).
 In particular, 
\[
\sp(n,\RR):=\Lie (Sp(n,\RR))=\{S:\ \trans SJ+JS=0\}.
\]
The Lie bracket in  $\sp(n,\RR)$ is just the commutator 
\[
[S_1,S_2]=S_1S_2-S_2S_1.
\]
\begin{itemize}
\item If $T\in Sp(n,\RR)$, we identify $T$ with the automorphism 
\[
T(z,u):=(Tz,u)
\]
of $\HH_n$, so that $Sp(n,\RR)$ embeds into the automorphism group 
 $\Aut(\HH_n)$ of $\HH_n.$

\item Further   automorphisms are the (anisotropic) {\it dilations} 
\[
\delta_r(z,u):=(rz,r^2u), \ r > 0,
\]
and the ``{\it Cartan involution}''
\[
\theta(x,y,u):=(x,-y,-u).
\]
\end{itemize}

\begin{prop}%1.1
$\Aut(\HH_n)$ is generated by $Sp(n,\RR)$, the dilations $\delta_r$, 
the inner automorphisms and $\theta$.
\end{prop}

\setcounter{subsection}{1}
\subsection{Integration on $\HH_n$}

The Lebesgue measure $dg:=dzdu$ is a {\it bi-invariant Haar measure}
 on $\HH_n$, i.e. 
\[
\int_{\HH_n} f(hg)\, dg=\int_{\HH_n}f(gh)\, dg=
\int_{\HH_n} f(g)\, dg\hskip1cm \forall h\in \HH_n.
\]
The {\it convolution} of two suitable functions (or distributions)
 $f_1,f_2$ on $\HH_n$ is defined by 
\begin{eqnarray*}
f_1\star f_2(g)&:=&\int_{\HH_n} f_1(h)f_2(h^{-1}g)\, dh\\
&=&  \int_{\HH_n} f_1(gh^{-1})f_2(h)\, dh.
\end{eqnarray*}
Define the {\it reflection at the origin} and the {\it involution}
 of $f$ by
\[
\check{f}(g):=f(g^{-1}) \mbox{ and }\  f^*(g):=\o{f(g^{-1})}, 
\mbox{ respectively.}
\]
Then, for $^\sharp =\check{}, \, ^*$, one has  $(f^\sharp)^\sharp=f$,
 and 
\[
(f_1\star f_2)^\sharp =f_2^\sharp\star f_1^\sharp, \hskip.5cm  
||f^\nat||_{L^1}=||f||_{L^1}.
\]
Notice that the {\it group algebra} $L^1(\HH_n,+,\star,^*)$ is a
 non-commutative involutive Banach algebra.

\vfill\pagebreak

\noi {\bf Remarks.} (a) Identify $z=(x,y)\in \RR^n\times \RR^n$ with 
$(z_1,\dots,z_n):= (x_1+iy_1,\dots,x_n+iy_n)\in 
\CC^n$, and call $f$ {\it polyradial}, if $f(z,u)=
\tilde f(|z_1|,\dots,|z_n|,u)$ for some function $\tilde f$ on 
$\RR_+^n \times \RR$. Under this identification of the underlying 
manifold of $\HH_n$ with $\CC^n$, the 
$n$-torus
 $\TT^n=\{(e^{i\ph _1}, \dots e^{i\ph_n}): \ph_i\in [0,2\pi[ \, \}$ 
acts by (symplectic) automorphisms
 $(z_1,\dots,z_n,u)\mapsto (e^{i\ph_1}z_1,\dots,e^{i\ph_n}z_n,u)$ 
on $\HH_n$, and $f$ is polyradial iff
 $f\circ \tau=f \ \forall \tau \in \TT^n$. Consequently, since 
\[
(f_1\star f_2)\circ \al =(f_1\circ \al)\star (f_2\circ \al)
\]
for every $f_1,f_2\in L^1(\HH_n)$ and $\al\in \Aut(\HH_n)$ with 
$\det D\al=1,$
$$\ L_{\rm pr }^1 (\HH_n):=\{f\in L^1(\HH_n): 
f \ \mbox{ is polyradial}\  \}$$
 forms a subalgebra of $L^1(\HH_n)$.
Even more is true:

\setcounter{subsection}{1}
\begin{prop}%1.2 
$L_{\rm pr}^1(\HH_n)$ is a commutative involutive Banach algebra.
\end{prop}

\noi {\bf Proof.} If $f\in L_{\rm pr}^1(\HH_n)$, then 
$\check{f}=f\circ \theta$. Hence, for $f_1,f_2\in
 L_{\rm pr}^1(\HH_n)$,
\begin{eqnarray*}
f_1\star f_2&=&(\check{f}_2\star \check{f}_1)\check{}=((f_2\circ 
\theta)\star (f_1\circ \theta))\check{}\\
&=&((f_2\star f_1)\circ \theta)\check{} 
=((f_2\star f_1)\check{} )\check{} =f_2\star f_1.
\end{eqnarray*}

\hfill Q.E.D.
\bigskip

\noi (b) If one replaces $\TT^n$ by the unitary group $U(n)$ in this 
discussion, one finds 
in a similar way that the {\it radial} $L^1$-functions $f$, i.e. 
functions which depend only on
 $|z|:=(|z_1|^2+\dots +|z_n|^2)^{1/2}$ and $u$, form a commutative 
subalgebra  $L_{\rm r}^1(\HH_n)$ of
 $L^1(\HH_n)$.

\bigskip

For $(z,u)\in \HH_n$, define the so-called {\it Koranyi-norm} by  

\begin{equation}%1.4
|(z,u)|:=(|z|^4+16u^2)^{1/4}=||z|^2\pm 4iu|^{1/2}.
\end{equation}
It has the following properties (Exercise):
\be
\item[(i)] $|\delta_r g|=r|g|\qquad \forall  g \in \HH_n,\ r > 0.$
\item[(ii)] $|g|=0 \iff g=0$.
\item [(iii)] $|g^{-1}|=|g|$.
\item[(iv)] $|g_1g_2|\le |g_1|+|g_2| \quad \forall g_1,g_2\in \HH_n$.
\ee
In particular, $|\cdot |$ is a so-called {\it homogeneous norm}, and 
$d_K(g_1,g_2):=|g_1^{-1}g_2|$ is a left-invariant metric on $\HH_n$.

\begin{remark}
$\HH_n$, endowed with the Koranyi-metric $d_K$ and the Haar measure, 
forms a space of homogeneous type in the sense of Coifman and Weiss.
\end{remark}

Denote by 
\[
B_r(g):=\{h\in \HH_n:|g^{-1}h|<r\}
\]
the open ball of radius $r>0$ centered at $g\in \HH_n$. Then, by 
left-invariance and (i),
$$
|B_r(g)|=|B_r(0)|=|\delta_r(B_1(0))|=r^Q|B_1(0)|,
$$
where
\[
Q=2n+2
\]
is the so-called {\it homogeneous dimension} of $\HH_n$.
\setcounter{subsection}{2}
\subsection{Left-invariant differential operators on $\HH_n$}

A linear operator $T:\cS(\HH_n)\to \cS'(\HH_n)$ is called {\it left}
 respectively {\it right - invariant}, if 
\[
T(\la_g\ph)=\la_g(T\ph)\quad \mbox{ respectively } \quad T(\vr_g\ph)
 =\vr_g(T\ph)
\]
for every $g\in G, \ \ph\in \cS$, where $\la$ and  $\vr$ denote the 
{\it left-regular} and  
{\it right-regular action} 
\[
(\la_g\ph)(h):=\ph(g^{-1} h),\quad (\vr_g\ph)(h):=\ph(hg).
\]
$T$ is called {\it homogeneous of degree } $\al \in \CC$, if 
\[
T(\ph\circ \delta _r)=r^\al (T\ph)\circ \delta_r\qquad 
\forall r > 0, \ \ph\in \cS.
\]
\bigskip

\noi{\bf The Lie algebra $\h_n$ of $\HH_n$}
\medskip

\noi Identify the tangent space $T_0 \HH_n$ with $\RR^{2n}\times \RR$.
 For $X\in T_0 \HH_n$, let $L_X$ denote 
the Lie-derivative 
\[
(L_X \ph)(g):=\frac d{dt}   \ph(g\cdot \gamma(t))|_{t=0},
\]
where $\gamma:[0,1] \to \HH_n$ is any smooth curve with $\gamma(0)
=0, \ \dot\gamma(0)=X$. 
Then $L_X$ is a left-invariant vector field on $\HH_n$, and the mapping 
$X\to L_X$ is bijective from 
$T_0\HH_n$ onto the space of all left-invariant real vector fields 
on $\HH_n$. 
In particular, the {\it Lie bracket} $[\ ,\ ]$ on $T_0\HH_n$ can be 
defined by 
\[
L_{[X,Y]}=[L_X,L_Y]:=L_XL_Y-L_YL_X.
\]
$T_0\HH_n$, endowed with $[\ ,\ ]$, forms the Lie algebra $\h_n$ of
 $\HH_n$. As usually, we shall 
henceforth identify $X\in \h_n$ with the corresponding Lie derivative $L_X$.

One computes easily that a basis of $\h_n$ is given by the vector fields
\begin{eqnarray}%1.5
X_j:=\frac{\de}{\de x_j}-\frac 1 2 y_j\frac{\de}{\de u},\quad \ Y_j:=\frac{\de}{\de y_j}+
\frac 1 2  x_j \frac{\de}{\de u},
\quad \ j=1,\dots,n,\ \mbox{ and} \quad  U:= \frac{\de}{\de u}.\\\nonumber
\end{eqnarray}

These satisfy the ``{\it Heisenberg commutation relations}''
\begin{eqnarray*}
[X_j,Y_k]&=&\delta  _{jk}\, U,\\
\lbrack X_j,X_k]&=&[Y_j,Y_k]=0,\\
\lbrack X_j,U]&=&[Y_j,U]=0.
\end{eqnarray*}

\noi {\bf Observe:} The $X_j,Y_j$ are homogeneous of degree 1, the 
``central derivative'' $U$ is
 homogenous of degree 2.

\medskip

\noi If $n=1$, we shall often write $X,Y$ in place of $X_1,Y_1$.

\medskip

\noi Notice that the exponential mapping $\exp:\h_n\to \HH_n$ is the 
identity mapping. 

Denote by $\u(\h_n)$ the associative algebra of all left-invariant 
differential operators
 on $\HH_n$. $\u(\h_n)$ can be identified with the {\it universal 
enveloping algebra} of $\h_n$. In
 particular, it is generated by the elements of $\h_n$.

\section{Local solvability }
\vskip-0.3cm\hskip1cm (see e.g. \cite{H1})
\vskip.5cm

\noi Let $P=P(x,D)=\sum\limits_{|\al|\le m} a_\al(x) D^\al$ be a 
linear PDO on $\RR^d$ of order $m$,
 where $D^\al=D_1^{\al_1} \dots D_d^{\al_d}$, $D_j=\frac 1{2\pi i} \ 
\frac \de{\de x_j}$. Denote by 
\[
P_m (x,\xi)=\sum_{|\al|=m} a_\al(x)\xi^\al,\quad (x,\xi)\in \RR^d\times \RR^d,
\]
its {\it principal symbol}. Assume that the coefficients $a_\al$ are smooth.

$P$ is said to be {\it locally solvable (l.s.)} at $x^0\in \RR^d$ if 
there exists an open neighborhood 
$\Om$ of $x^0$, such that for every $f\in C_0^\infty(\Om)$ there exists
a distribution $u\in \D'(\Om)$ 
solving the equation
\setcounter{equation}{0}

\begin{equation}
Pu=f\qquad \mbox{ in } \Om.
\end{equation}
We call $P$ locally solvable (in $\RR^d$), if it is locally solvable at 
every $x^0\in \RR^d$.

\begin{remark} By the theorem of Malgrange/Ehrenpreis, every constant 
coefficient PDO is locally solvable.
\end{remark}

\begin{example} {\rm Consider the left-invariant complex vector field
\[
Z=X+iY \qquad \mbox{ on } \HH_1.
\]
This  is just the famous  {\it Lewy-operator}, historically the first 
example of a linear PDO which is nowhere locally solvable.}
\end{example}

\noi {\bf Observe:} A left-invariant PDO on a Lie group is l.s. at one
 point of the group iff it is l.s. 
at every other point.

\medskip

Shortly after Lewy's example, H\"ormander produced the following
\medskip

\begin{theorem}[H\"ormander's criterion]
Assume there exists $\xi^0\in \RR^d$ s.t. 
\medskip

 $$  (H) \hskip3cm  P_m(x^0,\xi^0)=0 \quad \mbox{ and }\quad 
\{\RE P_m,\IM P_m\}(x^0,\xi^0)\ne 0,
$$

\noi where 
\[
\{a,b\}:=\sum_{j=1}^d\left( \frac{\de a}{\de \xi_j}  \frac{\de b}
{\de x_j}-\frac{\de a}{\de x_j}  \frac{\de b}{\de \xi_j}\right)
\]
denotes the Poisson bracket of $a$ and $b$. Then $P(x,D)$ is not locally
 solvable at $x^0$.
\end{theorem}

\noi Recall that $\xi^0$ is called {\it characteristic } for $P$ at $x^0$, if
 $P_m(x^0,\xi^0)=0$. 
\medskip

The lengthy proof makes use of the following 

\begin{basic lemma}
The equation (2.1) can be solved in $\Om$ if and only if the following 
holds true:

For every relatively compact open subset $\Lambda \subset \Om$ (shortly: $\Lambda
 \Subset \Om$) there exist constants  $C$ and $k\in \NN$, s.t. for every 
$f,v\in C_0^\infty(\Lambda)$,
\begin{equation}
|\int f v \, dx|\le C\sum_{|\al|\le k} ||D^\al f||_2\ \sum_{|\beta| \le k} 
||D^\beta \ \trans P v||_2
\end{equation}
Here, $\trans  P$ denotes the {\it formal transposed} of $P$, defined by 
\[
\int v (Pu)\, dx=\int (\trans  Pv)u\,dx.
\]
\end{basic lemma}

\bigskip

\noi{\bf Proof.} The sufficiency of (2.2) follows  by Hahn-Banach (exercise).

\medskip
Conversely, if $Pu=f$ can be solved for every $f\in \D(\Om)$ by some 
$u\in \D'(\Om)$, then 
\medskip

\noi ($*$) $\hskip3cm \lan f,v\ran =\int fvdx=\lan u, \trans Pv\ran\quad \forall 
v\in \D(\Lambda)$.

\medskip

Consider $\lan f,v\ran$ as a bilinear form on $C_0^\infty (\o \Lambda)\times
 C_0^\infty (\Lambda)$, where $C_0^\infty (\o \Lambda)$ is a Frechet space with the 
topology induced by the semi-norms $||D^\al f||_ 2$, and where $C_0^\infty(\Lambda)$ 
is endowed with the metrizable topology induced by the semi-norms 
$||{D^\beta}\,  \trans Pv||_2$.

Obviously, $f\mapsto \lan f,v\ran$ is continuous for fixed $v$. 

The continuity of $v\mapsto \lan f,v\ran$, for fixed $f$, follows on the other 
hand by ($*$).

Thus, $(f,v)\mapsto \lan f,v\ran$ is separately continuous, hence continuous,
 by Banach-Steinhaus. This proves (2.2).

\hfill Q.E.D.
\bigskip

\begin{remark}
Condition (2.2) is equivalent to 
\begin{equation}
 ||v||_{(-k)} \le C||\trans P v||_{(k)},
\end{equation}
where $||f||_{(\al)} =(\int (1+|\xi|^2)^\al |\hat f(\xi)|^2d\xi)^{1/2}$ denotes 
the Sobolev-norm of order $\al$.
\end{remark}
\bigskip

\noi{\bf Illustration of the proof of Theorem 2.3 in the case of Lewy's 
operator $Z$}
\medskip

\noi Assume w.r. that $x^0=0$.

A first important step is to find, for a given characteristic $\xi^0$ at 0
 satisfying condition $(H)$, a complex phase function of the form 
\medskip

\begin{equation}
 w(x)=\xi^0 \cdot x + i\, \trans  x\cdot A\cdot x+  O (|x|^3),
\end{equation}

\noi where  $\RE A$ is a positive-definite matrix, such that, if possible,
\begin{equation}
\trans  P(x,D)e^{2\pi i\la w}=0\qquad \forall \la \gg 1.
\end{equation}

\noi (This cannot always be achieved in the strict sense, only asymptotically as 
$\la\to\infty$, 
 but a necessary condition is that $w$ satisfies the {\it ``eikonal equation''}
\[
P_m(x,\nabla w)=0.)
\]
If $P=Z$ is Lewy's operator, then one computes that the characteristic 
points at 0 are $(0,0,\mu^0)$, which satisfy (H) if and only if $\mu^0\ne 0$.

A  suitable  phase can here be constructed directly by means of the 
following observation: Let

\begin{equation}
 q_\pm (z,u):=|z|^2\pm 4iu
\end{equation}
be the expression appearing implicitely in (1.4). Then one computes that 
\begin{equation}
Zq_+=0,
\end{equation}
so that $Z(f\circ q_+)=0$ for every holomorphic function $f$. Since $\trans  Z=-Z$, we may 
thus choose $w$ such that
\[
2\pi i w=-q_+ +q_+^2=-4iu -(|z|^2+16u^2)+ O((|z|+|u|)^3)
\]
in (2.5), with $\mu^0=-2/\pi$.

Given this phase, put
\[
v_\la:=e^{2\pi i\la w}\chi,\quad f_\la:=\la^3\chi(\la\cdot),
\]
where $\chi\in \D(\HH_1)$ is supported where $ |z|+|u| < 2\eps$, and  $\chi\equiv 1$ 
in $|z|+|u|\le \eps$. Then, as $\la\to +\infty$, 
\begin{eqnarray*}
\int_{\HH_1}f_\la v_\la \, dg=\int\int \chi(z,u)\chi(z/\la,u/\la)e^{2\pi i\la w( z/\la , u/ \la )}
 \, dzdu\to \int\int \chi (z,u)e^{-4iu}dzdu=\hat \chi(0,-\mu^0).
\end{eqnarray*}
On the other hand, 
\[
\trans Zv_\la =e^{2\pi i\la w }\, \trans Z\chi,
\]
where $\trans Z\chi$ is supported in the region where 
$|z|+|u|\sim \eps$. If $\eps$  is sufficiently small, then, by (2.4), 
$\IM\, w\sim \eps^2$ in this region, hence $|e^{2\pi i\la w}|\sim e^{-\delta \la}$, 
for some $\delta > 0$. This easily implies 
\[
||f_\la||_{(k)} \cdot || \trans Zv_\la||_{(k)}\to 0 \ \mbox{ as }\  \la \to +\infty.
\]
Thus, if we choose $\chi$ s.t. $\hat \chi (0,-\mu^0)\ne 0$, we obtain a 
contradiction to (2.2). 

\hfill  Q.E.D.
\medskip

\noi Remark: In general, (2.4) cannot be satisfied exactly, and the proof 
becomes considerably more involved. 

\medskip
For homogeneous left-invariant PDO's on $\HH_n$, the following necessary criterion 
for local solvability has 
proven extremely useful (analogues hold on general homogeneous groups).
\medskip

\noi{\bf Theorem 2.6} \cite{CR}, \cite{M1}. {\it Let $P\in \u(\h_n)$ be 
homogeneous. If $P$ is locally solvable, then there exist a Sobolev-norm 
$||\cdot ||_{(k)}$ and a continuous ``Schwartz-norm'' $||\cdot ||_\cS$ on
 $\cS(\HH_n)$, s.t.
\begin{eqnarray}
|f(0)| \le ||f||_\cS^{1/2} \, ||\trans Pf||_{(k)}^{1/2} \quad \forall  
f\in \cS(\HH_n).
\end{eqnarray}}

\noi{\bf Corollary 2.7} \cite{CR} {\it  Suppose there exists a non-trivial 
$f\in \cS(\HH_n)$ s.t. 

\medskip
\noi (CR) $\hskip 3cm \trans  Pf=0$.
\medskip

\noi Then $P$ is not locally solvable.}

\bigskip

\noi{\bf Proof.} Let $Q$ be an elliptic, right-invariant Laplacian on $\HH_n$, 
and let $\Om$ be an open neighborhood of $0,m\ge 1$.Then, for $\ph \in \D(\Om)$,
 by Poincar\' e's inequality and standard elliptic regularity theory,
\[
|\ph (0)|\le C'||Q^m\ph||_{2} \le C ||Q^{m+k/2} \ph||_{(-k)},
\]
provided $\Om$ is chosen sufficiently small.
We choose $k$ is as in (2.3), and assume $k$  to be even.  Since $Q^{m+k/2}$ commutes 
with the left-invariant operator $ \trans P$, by (2.3) we have 
\begin{eqnarray*}
||Q^{m+k/2} \ph||_{(-k)} &\le& C||{Q^{m+k/2}}\,  \trans P\ph||_{(k)}\\
&\le&C'||\trans P\ph||_{(2m+2k)},
\end{eqnarray*}
i.e. there exists a $K\in \NN, \ C\ge 0$, s.t.
\begin{eqnarray}
|\ph (0)|\le C \, ||\trans P \ph||_{(K)} \quad \forall \ph \in \D(\Om).
\end{eqnarray}

\noi Rescaling, we may assume w.r. that $\Om=B_2$, where $B_r:=B_r(0)$. Let 
$\trans ¿ P$ be homogeneous of degree $q$. Choose $\chi\in \D(B_2)$ s.t. 
$\chi\equiv 1$ on $B_1$. Then, for $f\in \cS$, by (2.9)
\medskip

\begin{equation}
|f(0)|\le C\, ||\trans P(\chi(f\circ \delta_r))||_{(K)} 
\quad \forall r > 0.
\end{equation}

\noi But:
\begin{eqnarray*}
\trans P(\chi(f\circ \delta_r))&=&\chi\, \trans P(f\circ \delta_r)+R(f\circ \delta_r)\\
&=& r^q\chi\,  (\trans Pf)\circ \delta_r+R(f\circ \delta_r),
\end{eqnarray*}
where $R=[\trans P,\chi]$ is a PDO whose coefficients are supported in 
$\{1\le |x|\le 2\}$. Thus, for $r\ge 1$, 
\begin{eqnarray*}
\lefteqn{||\trans P(\chi(f\circ \delta_r))||_{(K)}}\\
&&\le Cr^A\{||\trans Pf||_{(K)} +\sum_{|\al|\le N} (\int\limits_{1 < |x|<2} |f^{(\al)}
 (\delta_rx)|^2dx)^{1/2}\},
\end{eqnarray*}
for some constants $A>0, N\ge 0$. Now, 
\begin{eqnarray*}
\int_{1 < |x|<2} |f^{(\al)} (\delta_r x)|^2dx &\le& r^{-B}\int_{1 < |x|<2} 
|\delta _r x|^B|f^{(\al)}(\delta_rx)|^2dx\\
&\le& r^{-B-Q}\int |x|^B|f^{(\al)}(x)|^2dx.
\end{eqnarray*}
Choosing $B$ s.t. $A-B-Q=-A$, we find a Schwartz-norm $||\cdot ||_\cS$ s.t.
$$
||\trans P( \chi(f\circ \delta_r))||_{(K)}\le C ( r^A 
||\trans  Pf||_{(K)} +r^{-A}||f||_\cS).
$$
Combining this with (2.10) and optimizing in $r$ we obtain (2.8) (if we assume w.r.
 that $|f(0)|\le ||f||_\cS)$.

\hfill Q.E.D.
\bigskip

\noi In order to apply the ``$(CR)$-test'' from Corollary 2.7, one needs to construct functions 
in the kernel of $\trans P$. Here, representation theory can help.

\section{The group Fourier transform }
\vskip-0.3cm\hskip1cm (see e.g. \cite{Fo}, \cite{CG}, \cite{T})
\vskip.5cm

\setcounter{equation}{0}
Let $G$ be a locally compact group and $\H$ a Hilbert space. A 
{\it unitary representation } of $G$ on $\H$ is a strongly continuous
 homomorphism 
\[
\pi:G\to U(\H)
\]
of $G$ into the group $U(\H)$ of unitary operators on $\H$. We shall also write 
$\H_\pi$ in place of $\H$, if we want to emphasize that $\H$ is the 
representation space of $\pi$. Two representations $\pi$ and $\rho$ are 
called {\it equivalent}, if there exists a linear isometry $T$ from $\H_\rho$  
onto $\H_\pi$ such that $\pi(g) T=T\rho(g)$ for every $g\in G.$ $\pi$ is 
called {\it irreducible}, if the only closed and $\pi(G)$-invariant 
subspaces of $\H$ are $\{0\}$ and $\H$. The unitary dual $\hat G$ of 
$G$ consists of all equivalence classes $[\pi]$ of irreducible unitary 
representations. Often one identifies $\hat G$ also with a system of 
representatives of representations. 

As a consequence of the Stone-von Neumann theorem, such a system is
 given for the Heisenberg group $\HH_n$ by the following 
irreducible representations:
\be
\item[(i)] For $\mu\in \RR^\times:=\RR\setminus\{0\}$, the 
{\it Schr\"odinger representation } $\pi_\mu$ acts on $L^2(\RR^n)$ 
as follows:
\begin{eqnarray}
[\pi_\mu(p,q,u)f](x):=e^{2\pi i\mu(u+q\cdot x +\frac 1 2 q\cdot p)}
 f(x+p), \quad f\in L^2(\RR^n).
\end{eqnarray}

\item[(ii)] For $\zeta\in \RR^{2n}$, the characters 
\[
\om_\zeta(z,u):=e^{2\pi i \zeta \cdot z}
\]
are  1-dimensional  representations of $\HH_n$.
\ee
The characters are the irreducible  representations which act trivially
 on the center
\[
Z_n:=\{(0,0,u):u\in \RR\}
\]
of $\HH_n$, and they will play no role in the discussions to follow.

If $\pi$ is a unitary representation of $G$, and if $f\in L^1(G,dg)$ 
($dg$= left-invariant Haar measure), one defines $\pi(f)\in \B(\H)$ by 
\[
\pi(f)\xi:=\int _G f(g)\pi(g)\xi \, dg, \qquad \xi\in \H.
\]
One checks that the operator norm of $\pi(f)$ satisfies 
$||\pi(f)||\le ||f||_{L^1}$, and that the following holds true:

\bigskip

\noi {\it The ``integrated'' representation $\pi$ is a continuous 
homomorphism 
\[
\pi:(L^1(G),+,\star,^*)\to (B(\H),+,\circ ,^*)
\]
of involutive Banach algebras.}
\bigskip

For $f\in L^1(G)$, we define the {\it (group-) Fourier transform }
 of $f$ as the mapping $\hat f :\hat G\to \dot{\bigcup}_{\pi\in \hat G}
 \B( \H_\pi)$, given by 
\[
\hat f(\pi):=\int f(g)\pi(g)^*\,  dg=\int f(g)\pi(g^{-1})\, dg.
\]
Observe that for instance for $G=\HH_n$, 
\[
\hat f(\pi)=\pi(\check f),
\]
which implies 
\begin{eqnarray}
(f_1\star f_2)^\wedge  (\pi)=\hat{f_2}(\pi)\circ \hat{f_1}(\pi).
\end{eqnarray}

On $\HH_n$, one has the following explicit {\it Fourier-inversion 
formula} for ''nice'' functions, such as for example  Schwartz-functions:
\begin{eqnarray}
f(g)=\int_{\RR^\times} \tr (\hat f(\pi_\mu)\pi_\mu(g))\ |\mu|^n\, d\mu,
\quad g\in \HH_n.
\end{eqnarray}
The corresponding {\it Plancherel-formula} reads as follows:
\begin{eqnarray}
\int_{\HH_n} |f(g)|^2\, dg=\int_{\RR^\times} ||\hat f(\pi_\mu)||_{\HS}^2 
|\mu|^n \, d\mu.
\end{eqnarray}

Here, $\tr A$ denotes the trace of the operator $A$, and $||A||_{\HS}:=(\tr A^* A)^{1/2}$ 
its Hilbert-Schmidt norm.

This holds for $f\in L^2(G)$ in a similar sense as in the 
Euclidean case. For $f\in L^1\cap L^2(G)$, where 
$\hat f(\pi_\mu)$ is well-defined for every $\mu\ne 0$, 
part of the statement is that $\hat f(\pi_\mu)$ is a 
Hilbert-Schmidt-operator for a.e. $\mu\in \RR^\times$.

Notice that the characters $\om_\zeta$ do not enter in these 
formulas.

Formulas (3.3) and (3.4) can be deduced from the Euclidean Fourier
 inversion formula as follows:

Direct computations, based on formula (3.1), show that
 $\hat f(\pi_\mu)$ can be represented as a kernel operator
\begin{eqnarray}
(\hat f(\pi_\mu)\ph)(x)=\int_{\RR^n} K_f^\mu (x,y)
\ph(y)\, dy,\qquad \ph\in L^2(\RR^n),
\end{eqnarray}
with integral kernel
\begin{eqnarray}
K_f^\mu (x,y)&=&\int\int f(x-y,q,u)e^{-2\pi i\mu(u+\frac q 2(x+y))} 
\, dqdu,\\\nonumber
&=& f(x-y,\widehat{\frac \mu 2(x+y)}, \hat \mu).
\end{eqnarray}
Since  $\tr \hat f(\pi _\mu)=\int K_f^\mu(x,x)\, dx$, (3.3) follows
 easily (Exercise).
\bigskip

\subsection*{The Fourier transform of a differential operator}

If $P\in \u(\h_n)$, then 
\[
P\ph=P(\ph\star\delta)=\ph\star(P\delta),\quad \ph\in \cS,
\]
i.e. $P$ can be represented by convolution from the right with the
 compactly supported distribution $P\delta$. But from (3.6), one 
sees that $K_f^\mu$ is well-defined as a tempered distribution kernel
 $K_f^\mu\in \cS'(\RR^n\times \RR^n)$ supported near the diagonal
 $x=y$, for every distribution $f\in \E'(\HH_n)$ with compact 
support. This implies that the integral operator (3.5), defined 
in the Schwartz-sense of distributions, is well-defined on 
$\cS(\RR^n)$, and 
\[
\hat f(\pi_\mu):\cS(\RR^n)\to \cS(\RR^n)
\]
is continuous for every $f\in \E'( \HH_n)$.

For $P\in \u(\h_n)$, we now define its Fourier transform by 
\[
\hat P(\pi_\mu):=\widehat{P\delta} (\pi_\mu):=\pi_\mu((P\delta)\check{} ).
\]
Approximating $P\delta$ by $P\delta  \star\ph_\eps \in \D$, where 
$\{\ph_\eps\}_{\eps > 0}$ denotes a Dirac sequence in $\D$, one 
finds from (3.2) that 
\begin{eqnarray}
\widehat{P\ph} (\pi_\mu)=\hat P(\pi_\mu)\circ \hat \ph 
(\pi_\mu),\qquad \ph\in \cS,
\end{eqnarray}
and 
\begin{eqnarray}
\widehat{AB} (\pi_\mu)=\hat A(\pi_\mu)\circ \hat B(\pi_\mu),
\qquad \forall A,B\in \u(\h_n),
\end{eqnarray}
since $(AB)\delta=A(B\delta \star \delta)=B\delta \star A\delta.$

Since $X_j\delta =\frac \de{\de x_j} \delta,\  Y_j\delta =
\frac \de {\de y_j}\delta,\ \ U\delta  =\frac \de {\de u}\delta$,\  we 
find from (3.6) that 
\begin{equation}
\hat X_j(\pi_\mu)=\frac \de {\de x_j},\quad \hat Y_j(\pi_\mu)
=2\pi i\mu x_j,\quad \hat U(\pi_\mu)=2\pi i\mu.
\end{equation}
Also, from (3.6), one sees that 
\begin{eqnarray}
K_{f\circ \delta_r}^{r^2\mu} (x,y)=r^{-n-2} K_f^\mu(rx,ry),\qquad r > 0.
\end{eqnarray}
If $P\in \u(\h_n)$ is homogeneous of degree $q$, then 
$f=P\delta$ satisfies $f\circ \delta_r=r^{-Q-q}f$, hence 
from (3.10) we get 
\begin{eqnarray}
K_{P\delta}^{r^2\mu} (x,y)=r^{q+n}K_{P\delta}^\mu (rx,ry).
\end{eqnarray}
From Corollary 2.7, we can now deduce 
\medskip

\begin{cor}
Let $P\in \u(\h_n)$ be homogeneous, and assume there exist
 $\mu^0\in \RR^\times$ and $\phi\in \cS(\RR^n),\ \phi\ne 0$, s.t. 
$\widehat{\trans P} (\pi_{\mu^0})\phi=0$. Then $P$ is not locally solvable.
\end{cor}

\noi{\bf Proof.} Assume for instance $\mu^0>0$. For $\mu >0$, put 
\[
\phi^\mu(x):=\phi\left(\left(\frac{\mu}{\mu^0}\right)^{1/2}x\right).
\]
Then, by (3.11), $\widehat{\trans P}(\pi_\mu)\phi^\mu=0 \qquad 
\forall \mu >0$.
Let $\chi\in C_0^\infty (\RR^+)$, and put 
\[
K^\mu(x,y):=\chi(\mu)\phi^\mu(x)\phi^\mu(y).
\]
From (3.6), it follows that $K^\mu=K_f^\mu$ for some unique 
function $f\in \cS(\HH_n)$. And, 
\[
(\trans Pf)^\wedge (\pi_\mu)=\widehat{\trans P}(\pi_\mu) 
\hat f(\pi_\mu)=0,
\]
since $\hat f(\pi_\mu)$ is represented by the kernel $K^\mu$. 
Thus, by Fourier inversion, $\trans Pf=0$.

\hfill Q.E.D.
\bigskip

\noi {\bf Example 3.2.}
By (3.9), for the Lewy operator $Z=X+iY$ on $\HH_1$, one has 
\[
\widehat {\trans Z}(\pi_{\mu^0})=-\left(\frac d{dx}+x\right),\quad\mbox{ if }
 \mu^0=-1/2 \pi.
\]
Thus, the Gaussian $e^{-x^2/2}$ lies in the kernel of 
$\widehat{\trans Z}(\pi_{\mu^0})$.

\bigskip
\noi{\bf Remark 3.3.}
For a representation theoretic pendant to Theorem  2.6, see \cite{M1}.

\bigskip

\section{ Second order PDO's on $\HH_n$ with real 
coefficients and the metaplectic group}
\setcounter{equation}{0}
In the remaining part of these lectures, we shall discuss 
the following (still largely open)
\medskip

\noi{\bf PROBLEM.} Classify all second order left-invariant PDO's 
on $\HH_n$ which are locally solvable.
\bigskip

Let me remark that local solvability has also been studied for 
operators of higher order, and on more general Lie groups, in
 particular for bi-invariant PDO's and for ``transversally elliptic''
 operators. Some reference to the vast  literature on the subject 
can be found in \cite{Ba} and \cite{MR5}.

We shall concentrate here on the case of homogeneous operators of 
degree 2, which are of the form 
\begin{eqnarray}
L=\sum_{j,k=1}^{2n} a_{jk}W_jW_k+i\al U,\qquad a_{jk}, \al \in \CC,
\end{eqnarray}
where $W_j:=X_j,\ W_{n+j} :=Y_j,\ j=1,\dots,n$.

Throughout this section, the $a_{jk}$ will be real; the case of 
complex coefficients will be discussed in the last section. For 
results in the non-homogeneous case, see e.g. \cite{MR4}, \cite{MPR2}.

Let us put $A:=(a_{jk})_{j,k=1,\dots,2n}$ and 
\begin{eqnarray}
S:=-AJ.
\end{eqnarray}
Observe that $A$ is real and symmetric   if and only if 
$S\in \sp (n,\RR)$. Since, as it turns out, solvability of 
the operator $L$ is very much ruled by the spectral properties 
of $S$, we shall put 

\[
\Delta_S:=\sum_{j,k=1} ^{2n} a_{jk}W_jW_k,\quad S\in \sp (n,\RR),
\]
where $A$ is related to $S$ by (4.2).

The following theorem gives a complete answer for operators of
 the form (4.1) and $A$ real (for a generalization to arbitrary 2-step
 nilpotent groups, see \cite{MR3}).
\bigskip

\begin{theorem}{\rm \cite{MR2}} The operator $L_\al:=\Delta_S+i\al U$ 
is not locally solvable if and only if all of the following three conditions 
 hold:
\be 
\item[(i)] $\al\in \RR$;
\item[(ii)] $S$ is semisimple and has purely imaginary spectrum 
$\sigma(S)$; in this case, there exists some $T\in \Sp(n,\RR)$ such
 that $S':=TST^{-1}$ takes on the normal form 
\begin{eqnarray}
S'=\left(
\begin{array}{ccc|ccc}
        &&&\la_1&&\\
        &0&&&\ddots &\\
        &&&&&\la_n\\
        \hline
        -\la_1&&&&&\\
        &\ddots & &&0&\\
        &&-\la_n &&&
\end{array}
\right),
\end{eqnarray}
with ''frequencies'' $\la_1,\dots,\la_n\in \RR.$
\item[(iii)] There are no constants $C,N >0$, s.t.
\begin{eqnarray} \left| \sum_{j=1}^n (2k_j+1)\la_j\pm \al\right| 
\ge C\ (1+k_1+\dots+k_n)^{-N}
\end{eqnarray}
for all $k_1,\dots,k_n\in \NN$.
\ee
\end{theorem}

\noi Before we discuss some of the methods employed in its proof, 
let us consider some examples:
\bigskip

\noi{\bf Example 1.} Assume $S$ is given by (4.3). Then 
\[
\Delta_S=-\sum_{j=1}^n \la_j(X_j^2+Y_j^2).
\]
If all $\la_j$ are of the same sign, $\Delta_S$ is a so-called 
{\it sub-Laplacian}. In this case, condition (4.4) is equivalent 
to $\al\not\in \C$, where $\C$ is the ``{\it critical set}'' 
\[
\C:=\{\pm \sum_{j=1}^n (2k_j+1)\la_j\ :\ k_j\in \NN\}.
\]
Observe that local non-solvability for these operators does not 
only depend on the principal part of order 2, but in fact in a 
crucial way on the first order part $i\al U$. This phenomenon, which 
is in sharp contrast to the behaviour of so-called ``principal type'' 
operators (see e.g.\cite{MR5}), had first been observed in the 
fundamental work \cite{FS} on the so-called Kohn-Laplacian 
$\Delta_K=\sum_{j=1}^n (X_j^2+Y_j^2)$. For general sub-Laplacians,
 see also \cite{BG}.

It is interesting to remark that the approach by Folland/Stein in [8]
 avoids representation theory. It is based on the explicit formula
\begin{eqnarray}
(\Delta _K+i\al U)\Phi_\al=\gamma _\al \delta,
\end{eqnarray}
where 
\[
\Phi_\al:=q_+^{-\frac{n-\al}{2}} q_-^{-\frac{n+\al}{2}},
\]
with $q_\pm(z,u)=|z|^2\pm 4iu$ given by (2.6), and 
\[
\gamma_\al:=\frac{c_n}{\Gamma\left(\frac{n+\al}{2}\right) 
\Gamma\left(\frac{n-\al} 2\right)}.
\]
Clearly, for $\al \not\in \C$, $\gamma_\al\ne 0$, hence 
$\frac 1{\gamma_\al} \Phi_\al$ is a fundamental solution of
 $\Delta_K+i\al U$, which implies its local solvability. 

This approach, however, is restricted to rather particular
 operators (compare also \cite{DPR}).
\bigskip

\noi{\bf Example 2.} $X_1^2+Y_1^2-\la(X_2^2+Y_2^2)$ on $\HH_2$ is 
locally solvable if and only if there are constants $C,N>0$ s.t.
\[
|\la-p/q|>C q^{-N}\qquad \mbox{ (compare (4.4))},
\]
for all odd $p,q\in \NN$, i.e. if and only if $\la$ is neither a 
rational number $p/q$ with odd $p$ and $q$, nor a Liouville number of ``odd type''.
\bigskip

\noi{\bf Example 3.} $X_1^2-Y_1^2+i\al U$ is locally solvable on 
$\HH_1$ for every $\al\in \CC$.

In fact, here $S=\pmatrix{0 & -1\cr -1 & 0}$, hence 
$\sigma(S)=\{-1,1\}$ is real.
\bigskip

\subsection*{Basic tools in the proof of Theorem 4.1}
\setcounter{subsection}{0}
\medskip

\subsection{``Symplectic'' changes of coordinates}
If $T\in Sp(n,\RR)\hookrightarrow  \Aut(\HH_n)$, then, since 
$\exp=\id$ for $\HH_n$, $X(f\circ T)(g)=\frac d{dt} f(T(g\exp tX)) |_{t=0}
=\frac d{dt}f(T(g)\exp tT(X))|_{t=0}=(T(X)f)(Tg)$ for every 
$X\in \h_n$. This implies (Exercise)
\begin{eqnarray}
\Delta_S(f\circ T)=(\Delta_ {TST^{-1}}f)\circ T.
\end{eqnarray}
Since $U(f\circ T)=(Uf)\circ T$, this shows that {\it solvability 
of $\Delta_S+i\al U$ depends only on the 
conjugacy class of $S\in \sp (n,\RR)$ under the real symplectic
 group $\Sp(n,\RR)$.}

\bigskip

\subsection{Application of the group Fourier transform}

Whereas H\"ormander's criterion cannot be used here to prove
 non-solvability, since $L_\al$ has a real principal symbol, 
Theorem 2.6 does apply in a very similar way as in Example 3.2\ .

Let us illustrate this in the case of the operators 
\begin{eqnarray}
L_\al=X^2+Y^2+i\al U \quad\mbox{ on } \HH_1.
\end{eqnarray}
By (3.9), we have
\[
\widehat {L_\al}(\pi_\mu)=\frac{d^2}{dx^2}-(2\pi\mu x)^2-2\pi\al\mu.
\]
But, $\frac{d^2}{dx^2}-(2\pi\mu x)^2$ is just a re-scaled 
{\it Hermite operator}, with eigenfunctions
\[
h_k^\mu(x):=(2\pi|\mu|)^{1/4}h_k((2\pi|\mu|)^{1/2}x)
\]
and associated eigenvalues
\[
-2\pi|\mu|\ (2k+1),\qquad k\in \NN.
\]
Here,
\[
h_k(x)=c_k(-1)^k e^{x^2/2}\frac{d^k}{dx^k}e^{-x^2}
\]
denotes the $L^2$-normalized Hermite function of order $k$.

Consequently, 
\begin{eqnarray}
\widehat {L_\al} (\pi_\mu)h_k^\mu=-2\pi|\mu|(2k+1+(\sign \mu)\al)h_k^\mu,
\end{eqnarray}
i.e. there exist $\mu$ and $k$ with $\widehat{ L_\al} (\pi_\mu)h_k^\mu=0$ iff
 $\al\in \C=\{\pm (2k+1):k\in \NN\}$.

So, by Corollary 3.1, $L_\al$ is not l.s., if $\al \in \C$.

In general, if $S$ satisfies (i), (ii) and (iii) in Theorem 4.1, then, 
by standard symplectic linear algebra, one finds $T\in \Sp(n,\RR)$, 
which conjugates $S$ into the form (4.3) (see e.g. \cite[Lemma 3.1]{MR2}, 
and assuming that $S=S'$, one can argue in a  similar way as above, 
keeping Remark 3.3 in mind.
\medskip

On the other hand, if conditions (i) and (ii) in Theorem 4.1 do apply, 
but the diophantine condition (4.4) fails, one can prove local 
solvability by means of the Fourier inversion formula (3.3).

On a formal level, and grossly oversimplifying compared to the general
 case, the argument, which we shall again demonstrate in the case of 
the operator (4.7), is as follows:

Suppose $\al\not\in \C$, i.e. that (4.4) fails. Then, by 
(4.8), the operator $\widehat {L_\al}(\pi_\mu)$ is invertible, with 
\begin{eqnarray}
||(\widehat {L_\al} (\pi_\mu))^{-1}||\le C|\mu|^{-1}.
\end{eqnarray}
Now, given $f\in \D(\HH_1)$, try to define a function $w$ on $\HH_1$ 
by putting 
\begin{eqnarray}
w(g):=\int_{\RR^\times}\tr (\widehat {L_\al}(\pi_\mu)^{-1}
\hat f(\pi_\mu)\pi_\mu(g))\ |\mu|\, d\mu.
\end{eqnarray}
Since then $\hat w(\pi_\mu)=\widehat {L_\al}(\pi_\mu)^{-1}\hat f(\pi_\mu)$, 
one finds that $(L_\al w)^\wedge(\pi_\mu)=\widehat {L_\al}(\pi_\mu)
\hat w(\pi_\mu)=\hat f(\pi_\mu)$, hence $L_\al w=f$ (at least on a formal level).

To make this argument rigorous, the main problem is that (4.10) will
 in general not converge, because of the blow-up of estimate (4.9) 
as $\mu\to 0$. This can be overcome as follows:

Define $v$ as $w$ by (4.10), only with $\widehat {L_\al}(\pi_\mu)^{-1}$ 
replaced by $2\pi i\mu \ \widehat {L_\al} (\pi_\mu)^{-1}$. Then $v$ turns 
out to be well-defined, and one finds that 
\begin{eqnarray}
L_\al v=Uf.
\end{eqnarray}
But, since $U$ is locally solvable, given any $\ph \in \D$, there 
is some $f\in \D$ s.t. $Uf=\ph $ on the support of $\ph $. But then 
\[
L_\al v=\ph\qquad\mbox{ on }\  \supp\, \ph,
\]
hence $L_\al$ is locally solvable.
\bigskip

\subsection{Twisted convolution and the metaplectic group}
\setcounter{subsection}{1}

For generic $S\in \Sp(n,\RR)$, the operator $\widehat {\Delta_S}(\pi_\mu)$ 
will no longer have a discrete spectrum, and the approach described
 above breaks down. 

What saves the day is the following

\begin{lemma}
For $S_1,S_2\in \sp (n,\RR)$, we have 
\[
[\Delta_{S_1},\Delta_{S_2}]=-2U\Delta_{[S_1,S_2]}.
\]
\end{lemma}

\noi {\bf Proof.} Exercise.
\bigskip

\noi Denote by $f^\mu$ the partial Fourier transform of
 $f$ ``along the center'' of $\HH_n$, i.e. 
\[
f^\mu(z):=\int_\RR f(z,u)e^{-2\pi i\mu u}du,\qquad \mu\in \RR.
\]
Moreover, for suitable functions or distributions $\ph,\psi$ 
on $\RR^{2n}$, define the $\mu$-{\it twisted convolution} of 
$\ph$ and $\psi$ by
\[
\ph \times_\mu\psi(z)=\int_{\RR^{2n}} \ph (z-z')\psi(z')
e^{\pi i\mu\lan  z-z',z'\ran}\,dz'.
\]
One easily verifies that, for suitable distributions 
$f_1,f_2$ on $\HH_n$,
\begin{eqnarray}
(f_1\star f_2)^\mu&=&f_1^\mu \times_\mu f_2^\mu,\nonumber\\
(f^*)^\mu&=&(f^\mu)^*.
\end{eqnarray}
One also easily sees that  $L^1(\RR^{2n},+,\times_\mu,*)$ is
 a (non-commutative) involutive Banach algebra, and (4.12) shows that 
$f\mapsto f^\mu$ is a $^*$-homomorphism of $L^1(\HH_n,+,\star ,^*)$ onto it
(another way to verifying these facts is by passage through the ''reduced'' 
Heisenberg group; compare \cite{Fo}).

\medskip
If $\mu=1$, we just speak of the {\it twisted convolution},
 and write $\ph \times \psi$ in place of $\ph \times_1\psi$.
\medskip

\begin{remark}
{\rm Twisted convolution shares many features of ordinary convolution.
 For example, one has Young's inequality
\[
||\ph  \times \psi||_{L^r}\le ||\ph||_{L^p}||\psi||_{L^q},
\]
if $1/p+1/q=1+1/r$. More surprising is the following fact 
(see \cite{Fo}): If $\ph,\psi\in L^2(\RR^{2n})$, then also 
$\ph\times\psi\in L^2(\RR^{2n})$, and 
\[
||\ph \times \psi||_{L^2} \le ||\ph||_{L^2}||\psi||_{L^2}.
\]
Now, if $P\in \u(\h_n)$, then from (4.12) we get  
\[
(Pf)^\mu=(f\star P\delta)^\mu=f^\mu \times_\mu(P\delta)^\mu,
\]
where clearly $(P\delta)^\mu$ is a distribution supported 
at $0\in \RR^{2n}$. This shows that there exists a PDO $P^\mu$ 
on $\RR^{2n}$ such that 
\begin{eqnarray}
(Pf)^\mu=P^\mu f^\mu,\quad f\in \cS(\RR^{2n}).
\end{eqnarray}
For instance, by (1.5),
\begin{equation}
X_j^\mu=\frac{\de}{\de x_j}-i\pi\mu y_j,\quad Y_j^\mu=\frac{\de}{\de y_j}+i\pi\mu x_j,\quad
U^\mu=2\pi i\mu.
\end{equation}
In particular, from Lemma 4.2, we get 
\[
[\Delta_{S_1}^\mu , \Delta_{S_2}^\mu]=-4\pi i\mu\, \Delta_{[S_1,S_2]}^\mu.
\]
Moreover, $\Delta_S^\mu$ is formally self-adjoint, hence the mapping 
\begin{eqnarray}
S\mapsto \frac i {4\pi \mu} \Delta_S^\mu
\end{eqnarray}
is a representation of $\sp (n,\RR)$ by (formally) skew-adjoint 
operators on $L^2(\RR^{2n})$.}
\end{remark}

Let us consider the case $\mu=1$. In \cite{Ho}, R. Howe has  
proved for this case that the map (4.15) can be exponentiated to a unitary 
representation of the 
 {\it metaplectic group} $M_p(n,\RR)$, a two-fold covering of 
the symplectic group. $M_p(n,\RR)$ can in fact be represented by 
twisted convolution operators of the form $f\mapsto f\times\gamma$, 
where the $\gamma$'s are suitable measures which, generically, are 
multiples of purely imaginary Gaussians 
\begin{eqnarray}
e_A(z):= e^{-i\pi\trans z\cdot A\cdot z},
\end{eqnarray}
with real, symmetric $2n\times 2n$ matrices $A$. In particular, one has 
\begin{eqnarray}
e^{i\frac t{4\pi }\Delta_S^1} f=f\times \gamma_{t,S},\quad t\in \RR.
\end{eqnarray}

The measures $\gamma_{t,S}$ have been determined explicitly in 
\cite{MR1}. To indicate how this can be accomplished, let us argue 
on a completely formal basis:

If $e_A,e_B$  are two Gaussians (4.16) such that 
$\det(A+B)\ne 0$, one computes that 
\[
e_A\times e_B=[\det (A+B)]^{-1/2} e_{A-(A-J/2)(A+B)^{-1}(A+J/2)},
\]
where a suitable determination of the root has to be chosen. 
Choosing $A=\frac 1 2 JS_1, \ B=\frac 1 2 JS_2$, with $S_1,S_2\in
 \sp (n,\RR)$, and assuming that $S_1$ and $S_2$ commute, one finds that

$$
e_{\frac 1 2 JS_1}\times e_{\frac 1 2 JS_2}=
2^n(\det (S_1+S_2))^{-1/2} \, e_{\frac 1 2 J[S_1S_2+I)(S_1+S_2)^{-1}]}.
$$
\medskip
This reminds of the addition law for the hyperbolic cotangent, namely 
\[
\coth (x+y)=\frac{\coth x\coth y+1}{\coth x+\coth y}.
\]
We are thus led to define, for non-singular $S$,
\begin{eqnarray}
A(t):=\frac 1 2 J\coth (tS/2),
\end{eqnarray}
 which is well-defined at least for $|t|>0$ small.

Then 
\[
e_{A(t_1)} \times e_{A(t_2)}=2^n
(\det (A(t_1)+A(t_2)))^{-1/2}\,  e_{A(t_1+t_2)}.
\]
And, from 
\[
\coth x+\coth y=\frac{\sinh (x+y)}{\sinh x\sinh y},
\]
we obtain (ignoring again the determination of roots)
\begin{eqnarray*}
[\det \sinh ((t_1+t_2)S/2)]^{1/2}\, [\det (A(t_1)+A(t_2))]^{-1/2}\\
= [\det \sinh (t_1S/2)]^{1/2}\, [\det \sinh (t_2S/2)]^{1/2}.
\end{eqnarray*}
Together this shows that
\begin{eqnarray}
\gamma_{t,S}:=2^{-n}[\det\sinh (tS/2)]^{-1/2}e_{A(t)}
\end{eqnarray}
forms a (local) 1-parameter group under twisted convolution,
 and it is not hard to check that its infinitesimal generator 
is $\frac i {4\pi} \Delta_S^1$.
\medskip

\noi{\bf Warning:} Formula (4.18) only holds true for 
``generic'' $S\in \sp (n,\RR)$ and $t\in \RR$.

\medskip
\noi If one defines the {\it symplectic Fourier transform } 
of $f$ on $\RR^{2n}$ by 
\[
\stackrel{\triangle}{f}(\zeta):=\int_{\RR^{2n}} 
f(z)e^{-i\pi\lan \zeta,z\ran} dz=\hat f(\frac 1 2 J\zeta),
\]
one obtains a formula analogous to (4.19) for 
$\stackrel{\triangle}{\gamma}_{t,S}$:
\begin{eqnarray}
\stackrel{\triangle}{\gamma}_{t,S}=
[\det \cosh (tS/2)]^{-1/2}e_{B(t)},
\end{eqnarray}
where
\begin{eqnarray}
B(t):=\frac 1 2 J\tanh (tS/2).
\end{eqnarray}

\subsection{Solvability of $L_\al$ if  $\sigma(S)\subset
 \CC\setminus (i\RR)$ and $\al \in \RR$.}

In this case, formulas (4.19), (4.20) do apply in the strict sense.
 Observing that $\lan \gamma_{t,S}, \ph\ran =2^{-2n}\lan \stackrel
{\triangle}{\gamma}_{t,S},\stackrel{\triangle}{\ph}\ran$, they 
easily imply that there are a Schwartz norm $||\cdot ||_\cS$ and
 a constant $\beta\ne 0$, s.t.
\begin{eqnarray}
|\lan \gamma_{t,S},\ph\ran|\le \frac 1{\cosh \beta t}||\ph||_\cS.
\end{eqnarray}
For arbitrary $\mu\ne 0$, put 
\begin{eqnarray}
\gamma_{t,S}^\mu (z):=\cases{\mu^n\gamma_{t,S}(\mu^{1/2}z), \quad \mu >0,&\cr
&\cr
        |\mu|^n\o{\gamma_{t,S}(|\mu|^{1/2}z)},\quad \mu < 0.}
\end{eqnarray}
Then one verifies (see \cite{MR1}) that 
\begin{eqnarray}
e^{i\frac{t}{4\pi \mu}\Delta_S^\mu}f=f\ \times_\mu\ \gamma_{t,S}^\mu.
\end{eqnarray}
Now, the idea to solve the equation
\begin{eqnarray}
L_\al F=(\Delta_S+i\al U)F=f
\end{eqnarray}
is as follows: By taking a partial Fourier transformation, (4.25) 
is equivalent to 
\[
\left(\frac i {4\pi \mu}\Delta_S^\mu-\frac{i\al}{2}\right) F^\mu 
=\frac i{4\pi\mu} f^\mu\quad \forall \mu\in \RR^\times.
\]
Formally, we then obtain $F^\mu$ by 
\[
F^\mu=-\int_ 0^\infty e^{\frac{it}{4\pi \mu} 
\Delta_S^\mu-\frac{i \al}2 t} \left( \frac i {4\pi \mu} f^\mu\right) dt,
\]
hence

$$
F(z,u)= -\int_{\RR^\times} \int _0^\infty e^{-i\frac \al 2 t}
 f^\mu \times_\mu \gamma_{t,S}^\mu(z)\,dt\, \frac{e^{2\pi i\mu u}}
{4\pi i\mu} \,d\mu= f\star K(z,u),
$$

where $K$ is the distribution, formally defined by 
\[
K(z,u)=-\int_{\RR^\times} \int_0^\infty e^{-i\frac \al 2 t} 
\gamma_{t,S}^\mu (z)\, dt \frac{e^{2\pi i\mu u}}{4\pi i \mu}\, d\mu.
\]
This suggests to \emph{define} $K$ by 
\begin{eqnarray}
\lan K,\ph\ran =-\int_0^\infty \int_{\RR^\times} \left\lan 
\gamma _{t,S}^\mu, \ph^{-\mu}\right\ran e^{-i\frac \al 2 t}
\, \frac{ d\mu}{4\pi i\mu}\,dt.
\end{eqnarray}

\noi Now from (4.22) one derives that there are constants
 $N,M\in \NN$ s.t.
\[
\left|\left\lan \gamma_{t,S}^\mu, \ph^{-\mu}\right\ran \right| 
\le C\frac{(1+|\mu|)^N}{(\cosh \beta t)|\mu|^M} ||\ph ^{-\mu}||_\cS.
\]
This estimate implies that the distribution $\tilde K$, defined 
in the same way as $K$, only with $d\mu$ replaced by 
$(2\pi i\mu)^{M+1}d\mu$, is in fact well-defined. Moreover, we 
then have
\[
L_\al (f\star \tilde K)=U^{M+1} f.
\]
From here on, one can argue in a similar way as in \S 4.2 to 
show that $L_\al$ is locally solvable. 
\bigskip

\noi{\bf Remark.} The discussion of the remaining cases in 
Theorem 4.1 requires considerably more care (see \cite{MR2}).

\section{Second order PDO's on $\HH_n$ with complex coefficients}
\setcounter{equation}{0}
The classification of locally solvable PDO's on $\HH_n$ of the form 
\begin{eqnarray}
L=\sum_{j,k=1}^{2n} a_{jk}W_jW_k +\mbox{ lower order terms }
\end{eqnarray}
with complex coefficients $a_{jk}$
appears to be a challenging problem, which as of yet has only been
 answered for a few classes of operators (see \cite{DPR}, \cite{MPR1},
 \cite{MPR2} and \cite{KM}). Let us briefly survey those results.

We write the principal part of $L$ again as $\Delta_S$, 
however, now with $S\in \sp (n,\CC)$, i.e. $S=S_1+i S_2$, with 
$S_1,S_2\in \sp (n,\RR)$. The operators studied in \cite{DPR}, 
\cite{MPR1}, \cite{MPR2} can be  described as follows:

Assume $\RR^{2n}$ decomposes into symplectic subspaces 
\begin{eqnarray}
\RR^{2n} =V_1 \oplus \cdots \oplus V_r,
\end{eqnarray}
where each $(V_j,\om_j)$, with $\om_j:= \lan \ ,\ 
\ran |_{V_j\times V_j}$, is a symplectic vector space, and 
where the $V_j$'s are pairwise orthogonal w.r. to $\lan \ ,\ \ran$. 

Moreover, assume that each $V_j$ is $S$-{\it invariant}, i.e. 
that $S_i(V_j)\subset V_j$ for $i=1,2$. Recall that a basis
 $e_1,\dots,e_m,\ f_1,\dots,f_m$ of a symplectic vector 
space $(V,\om)$ is
 called {\it canonical} or {\it symplectic}, if 
\[
\om (e_j,e_k)=\om(f_j,f_k)=0,\quad \om(e_j,f_k)=\delta_{jk}.
\]
Then, choosing such a basis for each subspace $V_j$, we assume 
that $S$ can be written as a block diagonal matrix 
\begin{eqnarray}
S=\pmatrix{\gamma_1S_{(1)} &&&\cr
        & \gamma_2S_{(2)}&&\cr
        &&\ddots &\cr
        &&& \gamma_rS_{(r)}},
\end{eqnarray}
with $\gamma_j\in \CC^\times$ and 
\begin{eqnarray}S_{(j)}^2 =-I,\qquad j=1,\dots,r.
\end{eqnarray}
Observe that (5.3) generalizes the case (ii), formula (4.3), 
in Theorem 4.1, which appears to be of particular interest, to the 
complex setting.

We may and shall assume that each of the symplectic subspaces $V_j$ 
in (5.2) is {\it minimal} in the sense that it does not contain any 
proper $S$-invariant symplectic subspace.
\bigskip

\begin{theorem} {\normalsize \cite{MPR2}} If at least one of the minimal 
subspace $V_j$ has dimension $>2$, then $\Delta_S+P$ is not locally 
solvable for all first order (not necessarily invariant) differential
 operators $P$ with smooth coefficients.
\end{theorem}

This result is proved by means of H\"ormander's criterion Theorem 2.3: If 
we put again $S=-AJ$, then it follows from (1.5) that the principal 
symbol of $\Delta_S$ is given by 
\begin{eqnarray}
\sigma_S((z,u),(\zeta,\mu)):=-\trans (\zeta-\pi\mu Jz)\cdot  A\cdot 
(\zeta-\pi \mu Jz).
\end{eqnarray}
And, a straight-forward computation yields (compare Lemma 4.2)
\begin{eqnarray}
\{\sigma_S,\sigma_{S'}\}=4\pi\mu\sigma_{[S,S']} \qquad 
\forall S,S'\in \sp (n,\CC).
\end{eqnarray}
Thus, H\"ormander's criterion, applied to $\Delta_S+P$, just
 reads as follows:
\medskip

There is some $\zeta\in \RR^{2n}$ such that 
\medskip

\noi{\rm (H$'$)}\hskip3cm\qquad $\trans\zeta A_1\zeta =\trans\zeta_{ A_2} \zeta=0 
\mbox{ and } \trans\zeta A_3\zeta\ne 0$, 
\medskip

\noi where 
\[
A_1:=S_1 J,\  \ A_2:=S_2 J \ \mbox{ and }\  A_3:=[S_1,S_2]J.
\]

\bigskip

\noi{\bf Open Problem.} Classify all $S=S_1+iS_2\in \sp (n,\CC)$ 
for which (H$'$) applies. 
\medskip

In general, this seems to be a hard ``semi-algebraic'' problem.
 The proof of Theorem 5.1 makes use of a classification of 
normal forms of matrices
$S\in \sp (n,\CC)$ satisfying $S^2=-I$, with respect to
 conjugation by real symplectic matrices $T\in \Sp(n,\RR)$.
 Such a classification has been given in \cite{MT}. There remains the 
\bigskip

\subsection*{The case where all of the ``blocks'' $\gamma_jS_{(j)}$ 
are of size $2\times 2$}

According to the classification of normal forms in \cite{MT}, 
the $S_{(j)}$ can then be assumed to be either of the form
\begin{eqnarray}
S_{(j)} =\pmatrix{ i\eps_j\la_j & \la_j^2-1\cr
        1 & -i\eps_j\la_j} \qquad \mbox{ ``Type 1''},
\end{eqnarray}
with $\la_j\in \{-1\} \cup [0,\infty[$ and 
\[
\eps_j=\cases{ 1, & if $|\la_j|\le 1$,\cr
        \pm 1, & if $\la_j >1$,}
\]
or of the form
\begin{eqnarray}
S_{(j)}=\pmatrix{ 0 & i\cr i  & 0} \qquad \mbox{ ``Type 3''}.
\end{eqnarray}
The corresponding operators $\Delta_{S_{(j)}}$ are given by 
\begin{eqnarray}
L_{\la_j,\eps_j}:=(1-\la_j^2)X_j^2+Y_j^2-i\eps_j\la_j(X_jY_j+Y_jX_j)
\qquad \mbox{ ``Type 1''} 
\end{eqnarray}
and 
\begin{eqnarray}
-i(X_j^2-Y_j^2)\qquad \mbox{ ``Type 3''}.
\end{eqnarray}
\bigskip

\noi{\bf The case $n=1$}
\medskip

\noi Let us briefly discuss operators 
\begin{eqnarray}
L_\la :=(1-\la^2)X^2+Y^2+i\la(XY+YX)
\end{eqnarray}
on $\HH_1$. We call $L_\la$ a {\it generalized sub-Laplacian}, if 
$0 \le \la < 1$, i.e. if $\RE L_\la$ is a sub-Laplacian. In the case 
$\la=1$, i.e. 
\begin{eqnarray}
L_1=Y^2+i(XY+YX),
\end{eqnarray}
we speak of a {\it degenerate generalized sub-Laplacian}. If $\la >1$, 
then $\RE L_\la$ is more of ``hyperbolic type''.
\bigskip

\noi For $\HH_1$, local solvability of left-invariant operators (5.1) 
can be discussed in a complete way (\cite{MPR2}). To indicate the flavour
 of these results, let me highlight a few examples:
\bigskip

\noi{\bf Example 1.} If $L_\la$ is a generalized sub-Laplacian, then 
$L_\al+i\al U$ is l.s. if and only if 
\[
\al \not\in \C:=\{\pm (2k+1):k\in \NN\}.
\]
This extends the result for sub-Laplacians. 
\bigskip

\noi{\bf Example 2.} If $L_1=Y^2+i(XY+YX)$ is a degenerate generalized 
sub-Laplacian, then $L_1+i\al U$ is locally solvable if and only if 
\[
\al\not\in \C^+:=\{(2k+1):k\in \NN\}.
\]
For instance, for $\al =-1$ and $\al=1$, respectively, putting $\tilde 
Z:=Y+2iX$, one has 
\begin{eqnarray*}
L_1-iU&=& Y\tilde Z,\\
L_1+iU&=&\tilde ZY.
\end{eqnarray*}
Since $\tilde Z$ is of ``Lewy-type'', hence non-solvable, clearly $L_1+iU$
 cannot be solvable. The fact that $Y\tilde Z$ is locally solvable is more
 of a surprise (see\cite{MPR1}).
\bigskip

\noi{\bf Example 3.} If $\la >1$, then $L_\la+P$ is not locally solvable
 for every $P\in \u(\h_1)$ of order 1. This result cannot be obtained 
from H\"ormander's criterion, since this fails to apply for arbitrary 
operators $\Delta_S$ on $\HH_1$ (Exercise). It is proved in \cite{MPR2} 
by means of a variant of Corollary 3.1, which applies even 
to non-homogeneous operators.
\bigskip

\noi{\bf The case $n\ge 2$}
\medskip

\noi In this case, ``most'' of the operators $\Delta_S+P$ are locally 
non-solvable, as can be shown, with some effort, by means of 
 H\"ormander's criterion. The ``exceptional'' 
operators $\Delta_S$, to which (H$'$) does not apply, are listed in 
\cite[\S 6.1]{MPR2}.
There are five such exceptional classes, of which I want to mention 
two here:
\be
\item[(i)] 
On $\HH_n, \ n\ge 2$, ``positive combinations of generalized 
sub-Laplacians and of degenerate generalized sub-Laplacians'', 
more precisely 
\[
\Delta_S=\sum_{j=1}^m \gamma_j[(1-\la_j^2)X_j^2+Y_j^2+
i\la_j(X_jY_j+Y_jX_j)] +i \sum_{j=m+1}^n \beta_j(X_j^2-Y_j^2),
\]
where $|\la_j|\le 1,\ \gamma_j\in \CC^\times$, $\beta_j >0$, and where 
all of the quadratic forms
\[
\RE (\gamma_j[(1-\la_j^2)\xi_j^2+\eta_j^2+2i\la_j\xi_i\eta_j])
\]
are positive-semidefinite.
\item[(ii)] On $\HH_2$, for $\la > 1$, 
\begin{eqnarray*}
\lefteqn{\Delta_S=(1-\la^2)X_1^2+Y_1^2+i\la(X_1Y_1+Y_1X_1)}\\
&&+(1-\la^2)X_1^2+Y_1^2-i\la(X_2Y_2+Y_2X_2)
\end{eqnarray*}
\ee
\bigskip

\noi{\bf ad (i).} Observe that here the matrix $A=SJ$ satisfies 
$\RE\, A\ge 0$. Defining $B(t)$ as in (4.21), one then finds that
 $\RE (iB(-it))$ is positive semidefinite for every $t\ge 0$, so 
that $\stackrel \triangle \gamma_{-it,S}$, defined by (4.20), still
 remains a ``good'' Gaussian. As has been shown in \cite{MPR2}, this
 can be used to treat these operators by means of suitable 
modifications of the approach outlined in \S 4.4. In particular, 
one finds that  $\Delta_SS+i\al U$ is locally solvable for
every $\al$ not contained in the exceptional set 
\[
E:=\{\pm \sum_{j=1}^n \gamma_j(2k_j+1): \ k_1,\dots,k_n\in \NN\},
\]
provided $m=n$. 

If $\RE A$ is positive definite, this result follows also from the general theory
 of ''transversally elliptic'' partial differential operators; see e.g. \cite{H2},
\cite{BGH}.
\bigskip

\noi{\bf Open Problem.} Will $\Delta_S+i \al U$ be locally solvable 
for generic $\al \in \CC$, if $S\in \sp (n,\CC)$ and $\RE (SJ)$ is
 positive semidefinite, but not definite?
%\medskip

One can show that H\"ormander's criterion fails in this case (Exercise).
\bigskip

\noi{\bf ad(ii).} As has been proved recently in \cite{KM}, the 
operator $\Delta_S+P$ is locally solvable for arbitrary left-invariant lower 
order terms $P$.

In fact, the symplectic change of basis 
\begin{eqnarray*}
\tilde X_1:=Y_1-\sqrt{\la^2-1} X_2,&& \tilde Y_1:=Y_2+\sqrt{\la^2-1}X_1,\\
\tilde X_2:=Y_2-\sqrt{\la^2-1} X_1,&& \tilde Y_2:=Y_1+\sqrt{\la^2-1}X_2,
\end{eqnarray*}
transforms $\Delta_S$ into the operator
\[
Q_\la:=\left(1+\frac \la{2\sqrt{\la^2-1}}\right) DE 
-\frac \la{2\sqrt{\la^2-1}}{\o D} \ {\o E},
\]
where 
\[
D:=\tilde X_1-i\tilde X_2,\quad  E:=\tilde Y_2+i\tilde Y_1.
\]
We may thus reduce ourselves to the study of operators in $\u(\h_n)$, whose leading terms 
are of the form 
\begin{eqnarray*}
L_\A:=\sum_{j,k} \al_{jk}X_jY_k,
\end{eqnarray*}
 where  $\A=(\al_{jk})_{jk}$ is a complex $2\times 2$-matrix. Now, 
\[
\widehat {L_\A}(\pi_\mu)=2\pi i\mu\sum_{j,k=1}^2 \al_{jk} \frac \de 
{\de x_j} \circ x_k
\]
is homogeneous of degree 0. The second important property of 
$Q_\la$ is that $\widehat {Q_\la}(\pi_\mu)$ is elliptic away  from 
the origin, since its principal symbol is given by
\[
-4\pi \mu\left[(1+\frac \la {2\sqrt{\la^2-1}} )
(x_2+i x _1)(\xi_1-i\xi_2)-
        \frac \la{2\sqrt{\la^2-1}} \o{(x_2+ix_2)(\xi_1-i\xi_2)}\right].
\]
It has been proved  in \cite{KM} that a left-invariant operator on $\HH_2$ with 
leading term  $L_\A$ 
is locally solvable,  whenever $\widehat {L_\A}(\pi_\mu)$ 
is elliptic away from 0 and $\det \A \neq 0.$

In higher dimensions, such an ellipticity property of the 
operators $L_\A$ can never hold, which seems to explain 
why the exceptional operators of type (ii) only arize on $\HH_2$.

\medskip
\noi Mathematisches Seminar,  C.A.-Universit\"at Kiel, Ludewig-Meyn-Str.4,
 D-24098 Kiel, Germany\\
e-mail: mueller@math.uni-kiel.de
 \end{document}